\documentclass{amsart}

\copyrightinfo{2007}{American Mathematical Society}

\newtheorem{theorem}{Theorem}[section]
\newtheorem{lemma}[theorem]{Lemma}
\newtheorem*{acknowledgement}{Acknowledgement}
\theoremstyle{definition}

\numberwithin{equation}{section}

\begin{document}

\title{Triangularization of a\ Jordan Algebra of Schatten Operators}

\author{Matthew Kennedy}

\address{Department of Pure Mathematics \\
University of Waterloo \\
Waterloo, Ontario \\
Canada N2L 3G1}
\email{m3kennedy@uwaterloo.ca}

\thanks{Research supported by NSERC}
\subjclass[2000]{Primary 47A15; Secondary 17C65}

\date{February 15, 2007}

\begin{abstract}
We show that a Jordan algebra of compact quasinilpotent operators which
contains a nonzero trace class operator has a common invariant subspace. As a
consequence of this result, we obtain that a Jordan algebra of quasinilpotent
Schatten operators is simultaneously triangularizable.
\end{abstract}

\maketitle

\section{Introduction}

Engel's Theorem that a Lie algebra of nilpotent operators on a finite
dimensional space is triangularizable is a fundamental result in the theory of
Lie algebras. The problem of whether Engel's Theorem could be extended to a
Lie algebra of bounded quasinilpotent operators on a Banach space was first
raised by Wojty\'{n}ski, who proved the extension in the case when the
operators belong to a Schatten $p$-class \cite{bbwoj1}. After the construction
of a quasinilpotent operator with no nontrivial invariant subspaces by Read
\cite{bbread}, the scope of Wojty\'{n}ski's problem was suitably refined to
the case when the operators in the Lie algebra are compact, but it took many
years before Shulman and Turovski\u{\i} \cite{bbshu1}, in a remarkable paper,
proved this extension.

Jacobson \cite{bbjac} proved a strengthened form of Engel's Theorem which
implied that a Jordan algebra of nilpotent operators on a finite dimensional
space is triangularizable. In analogy with Engel's Theorem, it is natural to
ask whether Jacobson's result extends to a Jordan algebra of compact
quasinilpotent operators on a Banach space. The main result of this paper is
that a Jordan algebra of quasinilpotent Schatten operators is triangularizable.

\section{Preliminaries}

Let $\mathcal{H}$ be a complex Hilbert space, and let $\mathcal{B}%
(\mathcal{H})$ be the algebra of all bounded linear operators on $\mathcal{H}$
with the induced operator norm $\left\Vert \cdot\right\Vert $. For $p<\infty$,
we denote the set of all Schatten $p$-class operators on $\mathcal{H}$ by
$\mathcal{C}_{p}$. Each $\mathcal{C}_{p}$ is an ideal in $\mathcal{B}%
(\mathcal{H})$. An operator is said to be a Schatten operator if it is
contained in some $\mathcal{C}_{p}$. A beautiful account of the Schatten
operators is given by Ringrose \cite{bbrin}.

The operators in $\mathcal{C}_{1}$ are the trace class operators; they form a
Banach algebra with respect to the trace norm $\left\Vert \cdot\right\Vert
_{tr}$. For $A\in\mathcal{C}_{1}$ and $B\in\mathcal{B}(\mathcal{H})$, the
inequalities $\left\Vert A\right\Vert \leq\left\Vert A\right\Vert _{tr}$ and
$\left\Vert AB\right\Vert ,\left\Vert BA\right\Vert \leq\left\Vert
B\right\Vert \left\Vert A\right\Vert _{tr}$ always hold. By $tr(A)$ we denote
the trace of $A$, and by $\sigma(B)$ we denote the spectrum of $B$.

A Jordan algebra of operators $\mathcal{J}$ is a subspace of $\mathcal{B}%
(\mathcal{H})$ which is closed under the anticommutator product $\{A,B\}=AB+BA
$, for $A,B\in\mathcal{J}$; we may equivalently define $\mathcal{J}$ to be a
subspace of $\mathcal{B}(\mathcal{H})$ which is closed under taking powers.

Historically, the theoretical development of Jordan algebras has closely
corresponded to the theoretical development of Lie algebras, a fact which has
been utilized to prove results in both areas. Indeed, in what follows we will
require the following definitions and facts about Lie algebras.

A Lie algebra of operators $\mathcal{L}$ is a subspace of $\mathcal{B}%
(\mathcal{H})$ which is closed under the commutator product $[A,B]=AB-BA$, for
$A,B\in\mathcal{L}$. A Lie ideal $\mathcal{I}$ of $\mathcal{L}$ is a Lie
subalgebra of $\mathcal{L}$ such that $[A,B]\in\mathcal{I}$ whenever
$A\in\mathcal{L}$ and $B\in\mathcal{I}$.

Let $\mathcal{L}$ be a Lie algebra of operators. An element $A\in\mathcal{L} $
induces a linear transformation $ad(A)$ on $\mathcal{B}(\mathcal{H})$ defined
by $ad(A)B=AB-BA$, for $B\in\mathcal{B}(\mathcal{H})$. Note that $\mathcal{L}$
is invariant under $ad(A)$; by $ad_{\mathcal{L}}(A)$ we denote the restriction
of $A$ to $\mathcal{L}$. The operator $ad_{\mathcal{L}}(A)$ is called the
adjoint representation of $A$. If $ad_{\mathcal{L}}(A)$ is quasinilpotent for
all $A\in\mathcal{L}$, then we say that $\mathcal{L}$ is an Engel Lie algebra.
A Lie ideal $\mathcal{I}$ of $\mathcal{{}L}$ is said to be an Engel ideal if
it is an Engel Lie algebra.

For a Jordan algebra of operators $\mathcal{J}$, the Lie algebra $\mathcal{L}
$ generated by$\mathcal{J}$ coincides with the linear span of $\mathcal{J}%
+[\mathcal{J},\mathcal{J}]$, where $[\mathcal{J},\mathcal{J}]$ denotes the
linear span of the set $\{[A,B]:A,B\in\mathcal{J}\}$. Using the identity%
\[
\lbrack A,[B,C]]=\{\{A,B\},C\}-\{\{A,C\},B\}\text{,}%
\]
it is easy to verify that $\mathcal{L}$ is closed under the commutator product.

\section{Triangularization of a Set of Operators}

A\ set of bounded linear operators on a Hilbert space is said to be reducible
if there is a non-trivial subspace invariant under all the operators in the
set. The set is said to be triangularizable if there exists a maximal subspace
chain consisting of subspaces which are invariant under all the operators in
the set. This definition of triangularization generalizes the familiar finite
dimensional notion of triangularization to infinite dimensions \cite{bbrad2}.

If $\mathcal{S}$ is a set of operators and $\mathcal{U}$ and $\mathcal{V}$ are
invariant subspaces for $\mathcal{S}$ with $\mathcal{V}\subseteq\mathcal{U}$,
then $\mathcal{S}$ induces a set $\mathcal{\bar{S}}$ of quotients with respect
to the quotient space $\mathcal{U}/\mathcal{V}$, where for $A\in\mathcal{S}$,
the operator $\bar{A}\in\mathcal{\bar{S}}$ is defined on $\mathcal{U}%
/\mathcal{V}$ by%
\[
\bar{A}(x+\mathcal{V})=Ax+\mathcal{V}\text{.}%
\]
Since $\mathcal{U}$ and $\mathcal{V}$ are invariant under $A$, $\bar{A}$ is a
well-defined bounded operator on $\mathcal{U}/\mathcal{V}$. A property of a
set of operators is said to be inherited by quotients if every collection of
quotients of a set satisfying the property also satisfies the property.

Although the notion of reducibility is obviously much weaker than the notion
of triangularizability, the following lemma often reduces the problem of
\ showing a set of operators is triangularizable to the problem of showing
that the set is reducible.

\begin{lemma}
[The Triangularization Lemma \cite{bbrad2}]\label{trilem}Let $\mathcal{P}$ be
a property of a set of operators that is inherited by quotients. If every set
of operators on a space of dimension greater than one which satisfies
$\mathcal{P}$ has a nontrivial invariant subspace, then every such set is triangularizable.
\end{lemma}

\section{The Main Result}

Let $\mathcal{X}$ be a complex Banach space. For $A\in\mathcal{B}%
(\mathcal{X})$ and $\lambda\in\mathbb{C}$, let%
\[
\mathcal{E}_{\lambda}(A)=\{x\in\mathcal{X}:\lim_{n\rightarrow\infty}\left\Vert
(A-\lambda)^{n}x\right\Vert ^{1/n}=0\}\text{.}%
\]
The set $\mathcal{E}_{\lambda}(A)$ is called an elementary spectral manifold.
It is known \cite{bbriesz} that $\mathcal{E}_{\lambda}(A)$ is closed and
nonzero if $\lambda$ is an isolated point in the spectrum of $A$, and that $A$
is quasinilpotent if $\mathcal{E}_{0}(A)=\mathcal{X}$. Wojty\'{n}ski
\cite{bbwoj2} established the following result in his proof that a Lie algebra
of compact operators which is closed (with respect to the induced operator
norm) is either an Engel Lie algebra, or contains a nonzero finite rank operator.

\begin{lemma}
[Wojty\'{n}ski \cite{bbwoj2}]\label{finrank}Let $A\in\mathcal{B}(\mathcal{X})$
be compact. Then for nonzero $\lambda\in\mathbb{C}$, the set $\mathcal{E}%
_{\lambda}(ad(A))$ consists of nilpotent finite rank operators.
\end{lemma}

We use Wojty\'{n}ski's Lemma to obtain a stronger result in the case that the
Lie algebra consists of trace class operators.

\begin{lemma}
Let $\mathcal{L}$ be a Lie algebra of trace class operators such that
$\mathcal{L}$ is closed with respect to the trace norm. Then either
$\mathcal{L}$ is an Engel Lie algebra, or $\mathcal{L}$ contains a finite rank
operator. (Note that $\mathcal{L}$ is not required to be closed with respect
to the induced operator norm.)
\end{lemma}

\begin{proof}
In what follows, let $\mathcal{\hat{L}}$ be the closed subspace of
$\mathcal{C}_{1}$ obtained by endowing the elements of $\mathcal{L}$ with the
trace norm. For $A\in\mathcal{L}$, let $\hat{A}$ be the corresponding element
of $\mathcal{\hat{L}}$; then $ad(\hat{A})$ denotes the adjoint operator on
$\mathcal{C}_{1}$ induced by $\hat{A}$, and $ad_{\mathcal{\hat{L}}}(\hat{A})$
denotes the restriction of $ad(\hat{A})$ to $\mathcal{\hat{L}}$. We first show
that $ad_{\mathcal{L}}(A)$ is quasinilpotent whenever $ad_{\mathcal{\hat{L}}%
}(\hat{A})$ is quasinilpotent, and that the spectrum of $ad_{\mathcal{\hat{L}%
}}(\hat{A})$ is countable.

Suppose $ad_{\mathcal{\hat{L}}}(\hat{A})$ is quasinilpotent; then
$\mathcal{\hat{L}}\subseteq\mathcal{E}_{0}(ad(\hat{A}))$. Since
\[
\left\Vert ad(A)^{n}B\right\Vert ^{1/n}\leq\left\Vert ad(A)^{n}B\right\Vert
_{tr}^{1/n}%
\]
for all $B$ in $\mathcal{L}$, we deduce that $\mathcal{L}\subseteq
\mathcal{E}_{0}(ad(A))$. Since $\mathcal{E}_{0}(ad(A))$ is closed, it contains
the closure of $\mathcal{L}$; it follows that $ad_{\mathcal{L}}(A)$ is quasinilpotent.

Consider the left multiplication operator $L(\hat{A})$ induced by $\hat{A}$ on
$\mathcal{C}_{1}$. For $\lambda\in\mathbb{C}$, if $(\hat{A}-\lambda)^{-1}$
exists, then $L((\hat{A}-\lambda)^{-1})$ is the inverse of $L(\hat{A}%
-\lambda)$. Hence $\sigma(L(\hat{A}))\subseteq\sigma(\hat{A})$. Similarly,
$\sigma(R(\hat{A}))$ $\subseteq\sigma(\hat{A})$, where $R(\hat{A})$ is the
right multiplication operator induced by $\hat{A}$ on $\mathcal{C}_{1}$. Since
$ad(\hat{A})=L(\hat{A})-R(\hat{A})$, $\sigma(ad(\hat{A}))\subseteq\sigma
(\hat{A})-\sigma(\hat{A})$. Then, since $\sigma(\hat{A})=\sigma(A)$, and since
$A$ is compact and has countable spectrum, it follows that $ad_{\mathcal{\hat
{L}}}(\hat{A})$ also has countable spectrum.

Let $A\in\mathcal{L}$ be such that $ad_{\mathcal{L}}(A)$ is not
quasinilpotent. Then from above, $ad_{\mathcal{\hat{L}}}(\hat{A})$ is also not
quasinilpotent, and has countable spectrum. It follows that $\sigma
(ad_{\mathcal{\hat{L}}}(\hat{A}))$ is a nonempty countable compact set which
contains a nonzero isolated point $\lambda$, and hence that $\mathcal{E}%
_{\lambda}(ad(\hat{A}))$ contains a nonzero element of $\mathcal{\hat{L}}$.
Since
\[
\left\Vert (ad(A)-\lambda)^{n}B\right\Vert ^{1/n}\leq\left\Vert (ad(A)-\lambda
)^{n}B\right\Vert _{tr}^{1/n}%
\]
for all $B\in\mathcal{L}$, we deduce that $\mathcal{E}_{\lambda}(ad(A))$
contains a nonzero element of $\mathcal{L}$. Since, by Lemma \ref{finrank},
every operator in $\mathcal{E}_{\lambda}(ad(A))$ is of finite rank, we
conclude that $\mathcal{L}$ contains a nonzero finite rank operator.
\end{proof}

\begin{lemma}
\label{traceab}Let $\mathcal{J}$ be a Jordan algebra of quasinilpotent
operators, let $\mathcal{L}$ be the Lie algebra generated by $\mathcal{J}$,
and let $\mathcal{I}$ be the Lie ideal of $\mathcal{L}$ generated by the trace
class operators in $\mathcal{J}$. Then $tr(AB)=0$ for all $A,B\in\mathcal{I}$.
\end{lemma}

\begin{proof}
Recall that the trace class operators form an ideal in $\mathcal{B}%
(\mathcal{H})$. Since every element in $\mathcal{L}$ may be written as a sum
of (associative) words in the elements of $\mathcal{J}$ of length two or less,
it suffices to show that when $A\in\mathcal{J}$ is a trace class operator,
every word in the elements of $\mathcal{J}$ of length four or less which
contains $A$ has trace zero. Note that since every operator in $\mathcal{J}$
is quasinilpotent, every trace class operator in $\mathcal{J}$ has trace zero,
and in particular, $tr(A)=0$.

Suppose $B,C,D\in\mathcal{J}$. Since $AB+BA\in\mathcal{J}$, $tr(AB)=-tr(BA) $.
By the cyclic property of the trace, $tr(AB)=tr(BA)$, and it follows that
$tr(AB)=0$.

It is less straightforward to handle words of length three. Iterating the
identity%
\[
2ABA=\{\{A,B\},A\}-\{A^{2},B\}
\]
implies that $ABC^{2}BA,CBA^{2}BC\in\mathcal{J}$, and the identity%
\[
2(ABC+CBA)=\{C,\{B,A\}\}+\{A,\{B,C\}\}-\{B,\{C,A\}\}
\]
implies that $ABC+CBA\in\mathcal{J}$. Then, since%
\[
(ABC-CBA)^{2}=(ABC+CBA)^{2}-2(ABC^{2}BA+CBA^{2}BC)\text{,}%
\]
we obtain that $(ABC-CBA)^{2}\in\mathcal{J}$ . This implies that
$(ABC-CBA)^{2}$ is quasinilpotent, and hence that $ABC-CBA$ is quasinilpotent,
meaning $tr(ABC)=tr(CBA)$. Since $ABC+CBA\in\mathcal{J}$, $tr(ABC)=-tr(CBA)$,
and it follows that $tr(ABC)=0$.

Finally, to handle words of length four consider the identity
\[
ABCD+BCDA=\{A,B\}CD+BC\{A,D\}-B\{A,C\}D\text{.}%
\]
Since $ABCD+BCDA$ is a sum of words in $\mathcal{J}$ of length three,
$tr(ABCD)=-tr(BCDA)$ from above. Since, by the cyclic property of the trace,
$tr(ABCD)=tr(BCDA)$, it follows that $ABCD$ has trace zero. Applying the
cyclic property of the trace one last time to the words $AB$, $ABC$, and
$ABCD$ gives the desired result.
\end{proof}

Shulman and Turovski\u{\i}, in addition to their extension of Engel's Theorem,
have obtained a number of results regarding the existence of invariant
subspaces for Lie algebras of operators. We require the following reducibility criterion.

\begin{theorem}
[Shulman-Turovski\u{\i} \cite{bbshu2}]\label{nilideal}Let $\mathcal{L}$ be a
Lie algebra of compact operators. If $\mathcal{L}$ has a nonzero Engel ideal,
then $\mathcal{L}$ is reducible.
\end{theorem}

For $A\in\mathcal{B}(\mathcal{H)}$, and an isolated point $\lambda\in
\sigma(A)$, let $P_{\lambda}(A)$ denote the Riesz projection of $A$
corresponding to $\lambda$. For an isolated point $\lambda\in\sigma(ad(A))$,
we similarly let $P_{\lambda}(ad(A))$ denote the Riesz projection of $ad(A)$
corresponding to $\lambda$. The following lemma is a special case of a
technical result of Shulman and Turovski\u{\i}.

\begin{lemma}
\label{adproj}Let $A\in\mathcal{B}(\mathcal{H})$ be a finite rank operator.
Then for $\lambda\in\sigma(ad(A))$,%
\[
P_{\lambda}(ad(A))=\sum_{\substack{\alpha,\beta\in\sigma(A) \\\alpha
-\beta=\lambda}}L(P_{\alpha}(A))R(P_{\beta}(A))\text{,}%
\]
where $L(P_{\alpha}(A))$ and $R(P_{\beta}(A))$ denote the left and right
multiplication operators induced by $P_{\alpha}(A)$ and $P_{\beta}(A)$ on
$\mathcal{B}(\mathcal{H})$ respectively.
\end{lemma}

For a Lie algebra $\mathcal{L}$, let $[\mathcal{L},\mathcal{L}]$ denote the
linear span of the set $\{[A,B]:A,B\in\mathcal{L}\}$. It is clear that
$[\mathcal{L},\mathcal{L}]$ is a Lie ideal of $\mathcal{L}$; it is often
called the derived algebra of $\mathcal{L}$. The proof of the following
theorem is based on an argument of Cartan \cite{bbbour}.

\begin{theorem}
\label{cartan}Let $\mathcal{L}$ be a Lie algebra of trace class operators such
that $tr(AB)=0$ for $A,B\in\mathcal{L}$. Then every finite rank operator in
$[\mathcal{L},\mathcal{L}]$ is nilpotent.
\end{theorem}

\begin{proof}
Let $A\in\lbrack\mathcal{L},\mathcal{L}]$ be a finite rank operator; then
$\sigma(A)$ is finite. Define $T$ by%
\[
T=%
{\textstyle\sum\limits_{\lambda\in\sigma(A)}}
\bar{\lambda}P_{\lambda}(A)\text{,}%
\]
where $\bar{\lambda}$ denotes the complex conjugate of $\lambda$. Then $T$ and
$A$ commute, $\sigma(T)=\{\bar{\lambda}:\lambda\in\sigma(A)\}$, and
$\sigma(TA)=\{\bar{\lambda}\lambda:\lambda\in\sigma(A)\}$. To show that $A$ is
nilpotent, it therefore suffices to show that $tr(TA)=0$.

It is clear that%
\[
ad(T)=%
{\textstyle\sum\limits_{\lambda\in\sigma(ad(T))}}
\lambda P_{\lambda}(ad(T))\text{.}%
\]
Since $P_{\lambda}(T)=P_{\bar{\lambda}}(A)$ for all $\lambda\in\sigma(A)$, and
since $\sigma(T)=\{\bar{\lambda}:\lambda\in\sigma(A)\}$, $P_{\lambda
}(ad(T))=P_{\bar{\lambda}}(ad(A))$ by Lemma \ref{adproj}, which gives%
\[
ad(T)=%
{\textstyle\sum\limits_{\lambda\in\sigma(ad(A))}}
\bar{\lambda}P_{\lambda}(ad(A))\text{.}%
\]
Since $\mathcal{L}$ is invariant under $ad(A)$, it is invariant under
$P_{\lambda}(ad(A))$, and hence under $ad(T)$, which means that $[T,B]\in
\mathcal{L}$ for all $B\in\mathcal{L}$.

Since $A$ belongs to $[\mathcal{L},\mathcal{L}]$, we may write $A$ as a sum of
elements of the form $[B,C]$, for some $B,C\in\mathcal{L}$. Since%
\[
tr(T[B,C])=tr(C[T,B])=0
\]
by hypothesis, we deduce that $tr(TA)=0$, and hence that $A$ is nilpotent.
\end{proof}

\begin{theorem}
\label{trclred}Let $\mathcal{J}$ be a Jordan algebra of compact quasinilpotent
operators. If $\mathcal{J}$ contains a nonzero trace class operator, then
$\mathcal{J}$ is reducible.
\end{theorem}

\begin{proof}
Let $T\in\mathcal{J}$ be a nonzero trace class operator, let $\mathcal{L}$ be
the Lie algebra generated by $\mathcal{J}$, and let $\mathcal{I}$ be the Lie
ideal in $\mathcal{L}$ generated by $T$. Then $\mathcal{I}$ is a Lie algebra
of trace class operators. Let $\mathcal{K}$ be the closure of $\mathcal{I}$
with respect to the trace norm, and let $\mathcal{L}^{\mathcal{\prime}}$ be
the Lie algebra generated by $\mathcal{L}$ and $\mathcal{K}$. If
$\mathcal{L}^{\prime}$ is reducible, then $\mathcal{L}$, and hence
$\mathcal{J}$, are also reducible. We proceed by showing the reducibility of
$\mathcal{L}^{\prime}$.

We claim that $\mathcal{K}$ is a Lie ideal of $\mathcal{L}^{\prime}$. To
verify this, it suffices to show that for all $A\in\mathcal{L}$ and
$B\in\mathcal{K}$, $[A,B]\in\mathcal{K}$. Let $(B_{n})$ be a sequence in
$\mathcal{I}$ which converges to $B$ with respect to the trace norm. We have%
\begin{align*}
\left\Vert \lbrack A,B_{n}]-[A,B]\right\Vert _{tr} &  =\left\Vert
(AB_{n}-B_{n}A)-(AB-BA)\right\Vert _{tr}\\
&  \leq\left\Vert A(B_{n}-B)\right\Vert _{tr}+\left\Vert (B_{n}-B)A\right\Vert
_{tr}\\
&  \leq2\left\Vert A\right\Vert \left\Vert B_{n}-B\right\Vert _{tr}\text{,}%
\end{align*}
by the inequality $\left\Vert A(B_{n}-B)\right\Vert _{tr},\left\Vert
(B_{n}-B)A\right\Vert _{tr}\leq\left\Vert A\right\Vert \left\Vert
B_{n}-B\right\Vert _{tr}$. It follows that $[A,B_{n}]$ converges to $[A,B]$
with respect to the trace norm, and hence that $[A,B]\in\mathcal{K}$.

If $\mathcal{K}$ is an Engel Lie algebra, then $\mathcal{L}^{\prime}$ contains
a nonzero Engel Lie ideal, and $\mathcal{L}^{\prime}$ is reducible by Lemma
\ref{nilideal}. We may therefore suppose that $\mathcal{K}$ is not an Engel
Lie algebra; then $\mathcal{K}$ contains a nonzero finite rank operator by
Lemma \ref{finrank}. Consider the Lie ideal $[\mathcal{K},\mathcal{K}]$ of
$\mathcal{K}$; it is also a Lie ideal of $\mathcal{L}^{\prime}$. Indeed, for
$A\in\mathcal{L}^{\prime}$, and $B,C\in\mathcal{K}$,
\[
\lbrack A,[B,C]]=-[B,[C,A]]-[C,[A,B]]
\]
by the Jacobi identity, and $[C,A],[A,B]\in\lbrack\mathcal{K},\mathcal{K}]$.
Letting $\mathcal{F}$ be the set of finite rank operators in $\mathcal{L}%
^{\prime}$, it is clear that $\mathcal{F\cap}[\mathcal{K},\mathcal{K}]$ is
also a Lie ideal of $\mathcal{L}^{\prime}$.

By Lemma \ref{traceab}, $tr(AB)=0$ for all $A,B\in\mathcal{I}$; by the
continuity of the trace with respect to the trace norm, it follows that
$tr(AB)=0$ for all $A,B\in\mathcal{K}$, and hence by Theorem \ref{cartan} that
$\mathcal{F\cap}[\mathcal{K},\mathcal{K}]$ is an Engel Lie ideal. If
$\mathcal{F\cap}[\mathcal{K},\mathcal{K}]\not =0$, then it is a nonzero Engel
Lie ideal of $\mathcal{L}^{\prime}$ . Otherwise, if $\mathcal{F\cap
}[\mathcal{K},\mathcal{K}]=0$, the nonzero Lie ideal $\mathcal{F\cap K}$ of
$\mathcal{L}^{\prime}$ is commutative, and hence is a an Engel\ Lie ideal.
Either way, $\mathcal{L}^{\prime}$ has a nonzero Engel Lie ideal; by Lemma
\ref{nilideal} we conclude that $\mathcal{L}^{\prime}$ is reducible.
\end{proof}

\begin{theorem}
Let $\mathcal{J}$ be a Jordan algebra of quasinilpotent Schatten operators.
Then $\mathcal{J}$ is triangularizable.
\end{theorem}

\begin{proof}
It is clear that the property of being a Jordan algebra of quasinilpotent
Schatten operators is inherited by quotients. By the Triangularization Lemma,
it therefore suffices to show that $\mathcal{J}$ is reducible. If
$\mathcal{J}=0$, then $\mathcal{J}$ is trivially reducible, and if
$\mathcal{J}$ contains a nonzero trace class operator, then $\mathcal{J}$ is
reducible by Theorem \ref{trclred}; hence we may suppose that $\mathcal{J}$ is
nonzero, and that $\mathcal{J}$ doesn't contain any nonzero trace class operators.

Let $\mathcal{L}$ be the Lie algebra generated by $\mathcal{J}$. Since
$\mathcal{J}$ is nonzero, for some $p\geq1$, $\mathcal{J}\cap\mathcal{C}_{p} $
is nonzero. For $n\geq1$, it follows inductively from the identity%
\[
\lbrack A_{1},[A_{2},A_{3}]]=\{\{A_{1},A_{2}\},A_{3}\}-\{\{A_{1},A_{3}%
\},A_{2}\}\text{,}%
\]
that every $(2n-1)$-commutator $[A_{2n-1},[A_{2n-2},...,[A_{2},A_{1}%
]]]\in\mathcal{J}$, for all $A_{1},...,A_{2n-1}\in\mathcal{J}$. It is well
known that every (associative)\ word in the elements of $\mathcal{C}_{p}$ of
length at least $p$ is a trace class operator; since $\mathcal{J}$ doesn't
contain any nonzero trace class operators, it follows that for $n>p$, every
commutator in the elements of $\mathcal{J}$ which contains at least $n$
elements of $\mathcal{J}\cap\mathcal{C}_{p}$ is equal to zero.

Let $\mathcal{I}$ be the Lie ideal of $\mathcal{L}$ generated by
$\mathcal{J}\cap\mathcal{C}_{p}$. Observe that every $n$-commutator in the
elements of $\mathcal{I}$ may be written as a sum of $k$-commutators in the
elements of $J $, where $k\geq n$, and that each such $k$-commutator contains
at least $n$ elements of $\mathcal{J}\cap\mathcal{C}_{p}$. From above, this
implies that for $n>p$, every $n$-commutator in the elements of $\mathcal{I}$
is equal to zero. In particular, $(ad_{\mathcal{I}}(A))^{n}=0$ for all
$A\in\mathcal{I}$, which implies that $\mathcal{I}$ is a nonzero Engel ideal
of $\mathcal{L}$. By Theorem \ref{nilideal}, we conclude that $\mathcal{L}$,
and hence $\mathcal{J}$, is reducible.
\end{proof}

\begin{acknowledgement}
The author is grateful to Heydar Radjavi for his invaluable guidance and
support, and to Victor Shulman and Yuri Turovski\u{\i} for many helpful
suggestions and comments.
\end{acknowledgement}

\end{document}